\documentclass[10pt,a4paper,twoside]{amsart} 
\usepackage{amsmath}
\usepackage{amssymb} \usepackage{latexsym}
\usepackage{amscd}
\usepackage{fancyhdr}
\usepackage{mathrsfs} 
\usepackage{texdraw}

\newtheorem{definition}{\bf Definition}[section]

\newtheorem{lemma}[definition]{\bf Lemma}
\newtheorem{theorem}[definition]{\bf Theorem}
\newtheorem{proposition}[definition]{\bf Proposition}

\newcounter{remark}
\newcounter{remark:int}

\newcommand{\hX}{\widetilde{X}}
\newcommand{\hD}{\widetilde{D}}
\newcommand{\wee}{{\scriptscriptstyle \vee}}
\newcommand{\Db}{\mathfrak{Db}}
\newcommand{\Dbrd}[1]{\mathfrak{Db}_{\widetilde{X}}^{{\rm rd}, #1}}
\newcommand{\Crdt}{{\mathcal C}^{\rm rd}_{\widetilde{X}}}
\newcommand{\Crdtp}[1]{{\mathcal C}^{{\rm rd}, #1}_{\widetilde{X}}}

\newcommand{\widebar}[1]{\overline{#1}} 

\newcommand{\R}{{\mathbb R}} \newcommand{\C}{{\mathbb C}}
\newcommand{\N}{{\mathbb N}} \newcommand{\Q}{{\mathbb Q}}
 
\newcommand{\Ctd}{{\mathcal C}_{\widetilde{X}, \widetilde{D}}}
\newcommand{\Z}{{\mathbb Z}} 
\newcommand{\Hy}{{\mathbb H}}
\newcommand{\Hom}{{\rm Hom}} 
\newcommand{\cHom}{\mathcal{H}om}

\newcommand{\Epq}[2]{{\Omega^{\infty, (#1,#2)}_{\Xan}}}

\newcommand{\Etop}[1]{{\Omega^{\infty,#1}_{\widetilde{X}}}}

\newcommand{\EYpq}[2]{{\Omega^{\infty, (#1,#2)}_{Y^{\rm an}}}}

\newcommand{\vi}{\varphi} \newcommand{\ve}{\varepsilon}
\newcommand{\bve}{{\mathcal E}} \newcommand{\vt}{\vartheta}

\newcommand{\lto}{\longrightarrow} 
\newcommand{\Crd}{{\mathcal C}^{\rm rd}} 
\newcommand{\Hrd}{H^{rd}} 
\newcommand{\Hdr}{H_{dR}} 

\newcommand{\M}{{\mathcal M}}
 \renewcommand{\H}{\mathcal{H}}
\renewcommand{\O}{{\mathcal O}} 
\newcommand{\OXY}{\O_{\widehat{X|Y}}} 
\newcommand{\A}{{\mathcal A}_{\widetilde{X}}} 
\newcommand{\AD}{{\mathcal A}_{\widetilde{X}}^{<D}} 
\newcommand{\fA}[1]{{\mathcal A}_{\widehat{\widetilde{X} | #1}}}
\newcommand{\dr}[2]{{\rm DR}_{\widetilde{#1}}^{#2}}
\newcommand{\DR}{{\rm DR}_{\Xan}} 
 
\newcommand{\DRD}{{\rm DR}^{<D}_{\widetilde{X}}}

\newcommand{\Amod}{{\mathcal A}_{\widetilde{X}}^{{\rm mod }D}} 
\newcommand{\DRmod}{{\rm DR}_{\widetilde{X}}^{{\rm mod }D}}
\newcommand{\DRmd}[1]{{\rm DR}_{\widetilde{X}}^{{\rm mod }D,#1}}

\newcommand{\DRmY}{{\rm DR}_{\widetilde{Y}}^{{\rm mod }Z}}
\newcommand{\DRdY}{{\rm DR}_{\widetilde{Y}}^{<Z}}

\newcommand{\bR}{\boldsymbol{R}}
\newcommand{\Ass}[2]{\mathcal{A}_{\hX,#1}^{#2}}
\newcommand{\Asn}[1]{\mathcal{A}_{\hX}^{#1}}
\newcommand{\As}[2]{{\mathcal A}_{\widetilde{#1}}^{#2}}

\newcommand{\AmdY}{{\mathcal A}_{\widetilde{Y}}^{{\rm mod }Z}}

\newcommand{\Lotimes}{\overset{\mathbb L}{\otimes}}
\newcommand{\Smod}{\mathcal{S}^{{\rm mod }D}}
\newcommand{\SD}{{}^\wee\!{\mathcal{S}}^{<D}}
\newcommand{\SDo}{{\mathcal{S}}^{<D}}
\newcommand{\dolb}{\widebar{\partial}}
\newcommand{\Emd}{{\mathcal P}^{{\rm mod} D}}

\newcommand{\bdot}{\raisebox{-.3ex}{\mbox{\boldmath $\cdot $}}}

\newcommand{\Rd}{{\mathcal R}d}

\newcommand{\Uan}{U^{\rm an}}
\newcommand{\Xan}{X^{\rm an}}

\newcommand{\cR}{\mathcal{R}}

\newcommand{\ord}{{\rm ord}}
\newcommand{\fra}{\mathfrak{a}} \newcommand{\frb}{\mathfrak{b}}
\newcommand{\LE}{\mathbb{E}}
\newcommand{\Irr}{{\rm Irr}}

\makeatletter \renewcommand{\section}{\@startsection{section}{1}{\z@}%
  {-\baselineskip}{1\baselineskip}{\large\bf}}
\renewcommand{\subsection}{\@startsection{subsection}{2}{\z@}%
  {-\baselineskip}{0.5\baselineskip}{\normalsize\bf}}
\renewcommand{\subsubsection}{\@startsection{subsubsection}{3}{\z@}%
  {-\baselineskip}{0.5\baselineskip}{\normalsize\bf}} \makeatother

\numberwithin{equation}{section}

\begin{document}

\author[M.Hien]{Marco Hien}
\address{Marco Hien\\ NWF I -- Mathematik\\ Universit{\"a}t Regensburg\\ 93040
  Regensburg\\ Germany}
\email{marco.hien@mathematik.uni-regensburg.de}

\title[Periods for flat algebraic connections]{Periods for flat
  algebraic connections}

\subjclass{14F40, 14F10} \keywords{Period integrals, meromorphic
  connections, algebraic de Rham cohomology}

\maketitle

\begin{abstract}
In previous work, \cite{hien}, we established a duality between the
algebraic de Rham cohomology of a flat algebraic connection on a
smooth quasi-projective surface over the complex numbers and the rapid
decay homology of the dual connection relying on a
conjecture by C.~Sabbah, which has been proved recently by
T.~Mochizuki for algebraic connections in any dimension. In the
present article, we verify that Mochizuki's results allow to
generalize these duality results to arbitrary dimensions
also. 
\end{abstract}

\section{Introduction}

In \cite{hien}, we proved a duality theorem for the
algebraic de Rham cohomology of a flat algebraic connection $(E,
\nabla) $ on a smooth quasi-projective surface over the complex
numbers assuming a conjecture of C.~Sabbah (\cite{sabbah1}, Conjecture
I.2.5.1) on the existence of the so-called good formal structure after
birational pull-back of meromorphic connections on surfaces. Lately, a proof
of Sabbah's Conjecture in the case of algebraic connections (which lie
in the main focus of our work) has been achieved by T.~Mochizuki (see
\cite{mochizuki}). Moreover, in another astonishing work
\cite{mochizukiBig}, T.~Mochizuki is able to prove a far reaching
higher-dimensional generalization of this result, namely the existence of a good lattice
again after a suitable birational pull-back (loc.cit., Theorem
19.5). As in the one- and two-dimensional case, the canonical next
step into structural examinations of flat meromorphic connections
consists in the question of lifting the good formal properties to
asymptotic analogues. This step is also carried out in Mochizuki's
paper \cite{mochizukiBig} and results in the construction of Stokes
structures for flat meromorphic connections in arbitrary
dimensions. With these tools in hand, we are now able to construct the
period pairing of a flat algebraic connection between the algebraic de
Rham cohomology and the rapid decay homology which we introduced in
\cite{hien} and prove the perfectness of this period pairing
unconditionally in all dimensions. As a consequence, assuming rational
structures of the given geometric data over a subfield $k \subset \C $
as well as of the local system associated to the flat connection over
another subfield $F \subset \C $, we deduce a well-defined notion of
the {\bf period determinant} arising from the comparison of these
rational structures by means of the perfect period pairing. This
determinant will be an element in the quotient $\C^\times / k^\times
F^\times $. One of the motivations to study this determinant lies in
the mysterious analogies between flat connections over the
complex numbers and $\ell $-adic sheaves on varieties
over finite fields, according to which the period determinant is
expected to behave as the $\ve $-factor in the latter theory. If one
restricts to the subcategory of regular singular connections,
corresponding to the tamely ramified sheaves over finite fields, the
period determinant and its analogies to the $\ve $-factor has been
extensively studied by T. Saito and T. Terasoma in
\cite{saitotera}. With the present work, we hope to contribute to a
generalization of these lines of thought to the irregular singular
case. 

Let us give a short summary of the contents of this paper. The main object under
consideration is the period pairing between the algebraic de Rham
cohomology of the given flat algebraic connection and the rapid decay
homology -- analytic in nature -- of the dual connection. We will
discuss this construction in section \ref{sec:periodintegrals}. The
main result Theorem \ref{thm:periodpair}, the perfectness of the
period pairing, relies on a local duality between the algebraic de
Rham complex and an analytic de Rham complex of the dual bundle with
asymptotically flat coefficients. This duality will be proved in
Theorem \ref{thm:locduality} and Proposition \ref{prop:CrdDRD}
compares the asymptotically flat de Rham cohomolgy with the rapid
decay homology creating the link between the local duality and the
period pairing.

\section{Formal completion and real oriented blow-up}

We consider the following geometric situation. Let $U $ be a smooth
quasi-projective variety over $\C $ and $X $ a smooth projective
variety such that $D:=X \smallsetminus U $ is a divisor with normal
crossings. Introducing real polar coordinates around each irreducible
component of $D $ leads to the real oriented blow-up $\pi:
\widetilde{X} \to \Xan $. Locally at some $p \in D $, we can choose
complex coordinates $x_1, \ldots, x_d $ such that $p=0 $ and $D=\{ x_1 \cdots
x_k=0 \} $. Then $\pi $ reads as
$$
\pi: \big( [0, \ve) \times S^1 \big)^k \times Y \to \Xan \ , \ \big( (
r_\nu, e^{i \vt_\nu} )_{\nu=1}^k,  y \big) \mapsto
\big( (r_\nu \cdot e^{i \vt_\nu})_{\nu=1}^k, y \big) \ ,
$$
where $Y $ is a small analytic neighbourhood of $0
\in \C^{d-k} $, $y=(x_{k+1}, \ldots, x_d) $ and $\ve>0 $ a small
positive number. Restricted to $\pi^{-1}( \Uan ) $ the real oriented
blow-up is a homeomorphism by which we read $\Uan $ as a subspace of
$\widetilde{X} $. We write $\widetilde{\jmath}: \Uan \hookrightarrow
\widetilde{X} $ for the inclusion.

On $\widetilde{X} $, we consider the following sheaves of functions
(cp. \cite{sabbah1}, II.1.). Firstly, $\widetilde{X} $
carries the structure of a real manifold with boundary. Additionally,
the logarithmic differential operators $\widebar{x}_\nu
\widebar{\partial}_{x_\nu} $ act on the sheaf
$C^\infty_{\widetilde{X}} $ of $C^\infty $ functions. Let $\A $ be the
sheaf of functions which in the local siuation as above is given by
  $$
  \A := \bigcap_{\nu=1}^k \ker (\widebar{x}_\nu
  \widebar{\partial}_{x_\nu})  \cap \bigcap_{\nu=k+1}^d \ker
  \widebar{\partial}_{x_\nu} \subset
  C^\infty_{\widetilde{X}} \ .
  $$
Local sections in $\A $ for some open $\Omega \subset \widetilde{X} $ are
differentiable functions on $\Omega $ which are holomorphic on $\Omega \cap \Uan
$ and admit an asymptotic development in the higher-dimensional
analogue of Poincar\'e's asymptotic developments due to Majima (cp.
  Proposition B.2.1 in \cite{sabbah1} and \cite{majima}).

Next, if ${\mathcal P}_{\widetilde{X}}^{<D} $ denotes the sheaf
  of $C^\infty $-functions on $\widetilde{X} $ which are flat on
  $\pi^{-1}(D) $, i.e. all of whose derivations vanish on $\pi^{-1}(D)
  $ (cp. \cite{mal3}), we consider
  $$
  \AD := \A \cap {\mathcal P}_{\widetilde{X}}^{<D} \ .
  $$
Similarly, the local sections of $\AD $ are holomorphic restricted to
$\Uan $ and with vanishing asymptotic development (cp. Proposition
II.1.1.11 in \cite{sabbah1}). More precisely, if $u \in \AD( \Omega ) $ 
then for any compact $K \subset \Omega $ and any $N \in \N^k $, the
function $u $ satisfies
\begin{equation}\label{eq:ad=rd}
|u(x)| \le C_{K,N} \cdot |x_1|^{N_1} \cdots |x_k|^{N_k}  \mbox{ for all
  } x \in K \smallsetminus \pi^{-1}(D) \ ,
\end{equation}
in terms of local coordinates as above such that locally $D= \{
x_1 \cdots x_k =0 \} $. 

Due to Lemme II.1.1.18 in \cite{sabbah1}, we have the following fine
resolution of $\AD $
\begin{equation} \label{eq:ADresol}
\AD \hookrightarrow \big( {\mathcal P}_{\widetilde{X}}^{<D}
\otimes_{\pi^{-1} C_X^\infty} \pi^{-1} \Epq{0}{\bdot}, \dolb \big) \ ,
\end{equation}
where $\Epq{p}{q} $ denotes the sheaf of
$C^\infty $-forms of degree $(p,q) $ on $X$.

Analogously, if $\Amod $ denotes the sheaf of functions on $\widetilde{X}
$ which are holomorphic on $\Uan $ and of moderate growth along
$\widetilde{D} $, we have a fine resolution
\begin{equation} \label{eq:Amodresol}
\Amod \hookrightarrow \big( \Emd \otimes_{\pi^{-1} C^\infty_X}
\pi^{-1} \Epq{0}{\bdot} , \dolb \big) \ ,
\end{equation}
with $\Emd $ being the sheaf of $C^\infty $-functions on
$\widetilde{X} $ with moderate growth at $\pi^{-1}(D) $.

Finally, let $\Dbrd{-s} $ denote the sheaf of {\em rapid decay
distributions} on $\widetilde{X} $, i.e. the sheaf whose local section
on small open $V \subset \widetilde{X} $ are distributions
$$
\vi \in \Db_{\widetilde{X}}^{-s}(V) :=
\Hom_{\rm cont}( \Gamma_c(V, \Etop{s}), \C)
$$
on the space $\Etop{s} $ of $C^\infty $ differential forms on
$\widetilde{X} $ of degree $s $ with compact support in $\widetilde{X} $
satisfying the following condition: choosing coordinates $x_1,
\ldots, x_n $ on $X $ such that $D \cap V=\{ x_1 \cdots x_k =0 \} $, 
we require that for any compact $K \subset V 
$ and any element $N \in \N^k $ there are $m \in \N $ and $C_{K,N} >0
$ such that for any test form $\eta $ with compact support in $K $ the
estimate
\begin{equation} \label{eq:rddistrib} |\vi(\eta)| \le C_{K,N} \sum_i
  \sup_{|\alpha| \le m} \sup_K \{ |x|^N |\partial^\alpha f_i| \}
\end{equation}
holds, where $\alpha $ runs over all multi-indices of degree less than
or equal to $m $ and $\partial^\alpha $ denotes the $\alpha $-fold
partial derivative of the coefficient functions $f_i $ of $\eta $ in
the chosen coordinates.

\section{Good formal lattices and decompositions {\rm (after
    T. Mochizuki)}} 
\label{sec:goodcompact}

\subsection{Deligne-Malgrange lattices}

Let there be given a flat meromorphic connection $\nabla $ on the
locally free $\O_{\Xan}(\ast D) $-module $E $ or rank $r $, where $D
\subset X $ is a divisor with normal crossings. After chioce of a local
trivialization $E \cong (\O_{\Xan}(\ast D))^r $, the connection reads as
$\nabla= d + A $ with the connection matrix $A \in M(r \times r,
\Omega^1_{\Xan}(\ast D)) $. A change $T \in {\rm GL}_r( \O_{\Xan}(\ast
D)) $ of the trivialization transforms the connection matrix due to
the formula 
$$
A' := T^{-1} \, dT + T^{-1} A T \ .
$$
The local classification of these connections, i.e. of the
connection matrices up to this transformation, is a difficult task in
general.

The first major subdivision of flat meromorphic connections lies in the
distinction between regular singular connections and irregular
singular ones. Let us recall this notion in the given geometric
situation, i.e. with $D $ being a normal crossing divisor. Then, a flat
meromorphic connection $(E, \nabla) $ is {\bf regular singular}, if
there is a trivialization $E \cong ( \O_X(\ast D))^r $ such that the
resulting connection matrix has logarithmic poles along $D
$ at most, i.e. $A $ can be written as
$
A= \sum_{i=1}^k A_i (x) \,  \operatorname{dlog} x_i + \sum_{j=k+1}^d
A_j(x) \, dx_j 
$
with holomorphic matrices $A_i(x) \in M(r \times r, \O_{\Xan}) $. The
structure of regular singular connections is well understood (see
\cite{deligne}), in particular we know that each regular singular
connection is a succesive extension of rank one connections. For the
latter, one finds a basis vector $e $ over $\O_{\Xan} $ such that
$$
\nabla e = e \otimes \big( \lambda_1 \operatorname{dlog} x_1 + \ldots +
\lambda_k \operatorname{dlog} x_k \big)
$$
for some $\lambda \in \C^k $. Such a connection will be denoted by
$x^\lambda $. 

Next, the most elementary irregular singular connections are rank one
connections which in a suitable basis vector $e $ read as $\nabla e =
e \otimes d \fra $ for some $\fra \in \O_{\Xan}( \ast D) $, which
up to isomorphism depends on $\fra \operatorname{mod} \O_{\Xan} $
only. Such a connection will be denoted by $e^{\fra} $. 

Now, let $Y \subset D $ be a stratum in the natural stratification of
the normal crossing divisor $D $. Passing from $\O_{\Xan} $ to the 
formal completion $\OXY $ of $\O_{\Xan} $ along $Y $ leads to the
problem of formal classification, to which a extensive answer is given
in the case $\dim X=1 $ by the Levelt-Turrittin Theorem (see
\cite{mal1}, Chapter III). In \cite{sabbah1}, C. Sabbah investigated
the two-dimensional situation leading to a precise conjecture as well
as partial results in this direction. Recently, T. Mochizuki was able
to give a proof Sabbah's Conjecture and a higher-dimensional
generalization using a different apporach by exmamining the so-called
good lattices which go back to Malgrange (see
\cite{malgrange_reseau}). We will now explain Mochizuki's result about
these formal properties of meromorphic connections. For more details,
we refer to \cite{mochizukiBig}, chapter 5. 

Let us choose local coordinates such that
$D=\{x_1 \cdots x_k = 0 \} $ with the irreducible components $D_i=\{
x_i=0 \} $ and consider the
local situation, i.e. we put $X= \Delta^k \times Y $ where $\Delta^k $
is s small poly-disc in $\C^k $ centered at the origin and $Y $ is a
small neighbourhood of the origin in the remaining
variables. Accordingly, we will denote the 
first $k $ variables by $z_1, \ldots, z_k $ and the others by $y_1,
\dots, y_{d-k} $. Let $\le $ be the partial order
on $\Z^k $ given by $m \le n $ if $m_i \le n_i $ for all $i $. Now,
any meromorphic function $f \in \O_{\Xan}(\ast D) $ admits a Laurent
expansion $f= \sum_{m \in \Z^k} f_m(y) z^m $ with holomorphic
functions $f_m \in \O_Y $. The {\em order of $f $} will be the minimum
$$
\ord(f) = \min \{ m \in \Z^k \mid f_m \neq 0 \} \ ,
$$
assuming that this minimum with respect to $\le $ exists. T.Mochizuki
defines the notion of a {\bf good set of irregular values on $(X,D) $}
to be a finite set $S \subset \O_{\Xan}( \ast D)/ \O_{\Xan} $, such
that
\begin{enumerate}
\item $\ord(\fra) $ and $\ord(\fra-\frb) $ exist for all $\fra \neq \frb $ in $S $,
\item the set $\{ \ord(\fra -\frb) \mid \fra, \frb \in S \} $ is
  totally ordered with respect to $\le $ on $\Z^k $,
\item the leading terms $\fra_{\ord(\fra)}(p) $ and $(\fra
  -\frb)_{\ord(\fra -\frb)}(p) $ of the Laurent expansions are
  non-vanishing for all $p \in Y $.
\end{enumerate}
For any subset $I \subset \{1, \ldots, k\} $, let $I^c $ be its
complement and furthermore $D_I := \bigcap_{i \in I} D_i $ and $D(I):=
\bigcup_{i \in I} D_i $. The completion of $X $ along $D_I $
resp. $D(I) $ will is denoted by $\widehat{D}_I $
resp. $\widehat{D}(I) $. For a given subset $S \subset \O_{\Xan}(\ast
D) / \O_{\Xan} $ let $S(I):=\{ \fra \, {\rm mod } \, \O_{\Xan}(\ast D(I^c))
\mid \fra \in S \} $. 

Let $\LE $ be a lattice in $E $, i.e. a
locally free $\O_{\Xan} $-module such that $\LE \otimes \O_{\Xan}
(\ast D) = E $. Then $\LE $ is called an {\bf unramified good lattice}
at $p \in D $ if there is a good set of irregular values $S $ as above
such that for each $I $ we have a formal decomposition
\begin{equation} \label{eq:formdecomp}
(\LE, \nabla)|_{\widehat{D}_I} = \bigoplus_{\fra \in S(I)} ( {}^I
\widehat{\LE_\fra}, {}^I \widehat{\nabla_\fra}) \ ,
\end{equation}
such that
$$
\big( \widehat{\nabla_\fra} - d \fra \big) ( {}^I \widehat{\LE_\fra} )
\subset {}^I \widehat{\LE_\fra} \otimes \big( \Omega^1_{\Xan} ( \log
D(I)) + \Omega^1_{\Xan} (\ast D(I^c)) \big) \ .
$$
The set $S $ of irregular values is uniquely determined by the given
connection and hence it is denoted by $S:= \Irr(\nabla) $.

The given connection $(E, \nabla) $ is said to admit a {\bf good
  lattice} at $p \in D $, if for some coordinate neighbourhood
$U= \Delta^d $ of $p $ as above, there exists a ramification map
$$
\rho_e:\Delta^d \to U \ , \ (t_1, \ldots, t_d) \mapsto (t_1^e, \ldots,
t_k^e, t_{k+1}, \ldots, t_d)
$$
for some $e \in \N_0 $ such that $\rho_e^\ast (E, \nabla) $ admits an
unramified good lattice $\widetilde{\LE} $. With these notions,
T.~Mochizuki proves the following far reaching result
\begin{theorem}[{\bf T. Mochizuki}, \cite{mochizukiBig}, Theorem 19.5]
    \label{thm:mochiformal}
Let $(E, \nabla) $ be a flat meromorphic connection on $(X,D) $. Then
there exists a regular birational map $\vi:X' \to X $ such that the
pull-back $\vi^\ast(E, \nabla) $ admits a good lattice at any
$p \in \vi^{-1}(D) $.
\end{theorem}
Actually, the theorem in \cite{mochizukiBig} is more
precise in the sense that due to Malgrange's work
(\cite{malgrange_reseau}) one knows that generically there is a
''canonical lattice'', which Mochizuki calls the Deligne-Malgrange
lattice. Mochizuki proves that its extension to the whole $X $ is
again locally free over $\O_{\Xan} $ and a good lattice.

\subsection{Good decomposition in multisectors}

If we apply Mochizuki's Theorem \ref{thm:mochiformal} to the given flat
meromorphic connection $(E, \nabla) $ on $(X,D) $, we know that the
pull-back with respect to some birational map and some finite
ramification admits a good lattice $\LE $. It will be important
for our purposes to see that the resulting formal decomposition can be
lifted to an asymptotic one on small multisectors. In dimension
two, this is the ingredient of \cite{sabbah1}, II.2, in higher
dimensions it can be derived from Mochizuki's approach of good
lattices and can be found in \cite{mochizukiBig} also, see the remark
after the theorem below. 

Let $\fA{Z} $ denote the formal completion of $\A $ along
$\pi^{-1}(Z) $ for a closed subset $Z \subset X $ and $T_Z $ the
natural morphism $\A \to \fA{Z} $. We will consider mainly the case
where $Z \subset D $ is a union of local irreducible components of $D
$. Note that $\fA{Z}|_{\pi^{-1}(Z)} = \pi^{-1}
\O_{\widehat{X|Z}} $ in such a case. The following generalization of the
one-dimensional Borel-Ritt theorem tells us that the sequence $0 \to \AD \to \A
\stackrel{T_D}{\to} \fA{D} \to 0 $ is exact (see \cite{sabbah1},
II.1.1.16). The formal decompositions \eqref{eq:formdecomp} given by
the existence of an unramified good lattice can indeed be lifted in
the follwing sense:

\begin{theorem}[{\bf T.~Mochizuki}, \cite{mochizukiBig}]
  \label{thm:mochiasymp} 
Let $(E, \nabla) $ be a flat meromorphic connection and assume that it
admits an unramified good lattice $\LE $. Let $\vt \in \pi^{-1}(D) $
be any multi-direction. Then there exists a small mutlisector $S
\subset X \smallsetminus D $ around $\vt $ with closure $\widebar{S}
\subset \widetilde{X} $ such that we have a $\nabla $-flat decomposition
\begin{equation} \label{eq:asymplift}
(\pi^\ast \LE \otimes_{\pi^\ast \O_{\Xan}} \A)|_{\widebar{S}} =
\bigoplus_{\fra \in \Irr(\nabla)} \mathcal{G}_{\fra, S} \ ,
\end{equation}
where the connection on $\mathcal{G}_{\fra,S} $ induced by $\nabla
- d \fra $ is logarithmic, hence regular singular.
\end{theorem}

The proof of Theorem \ref{thm:mochiasymp} is completely contained in
\cite{mochizukiBig}. However, since Mochizuki develops the theory of
Stokes structures in a more general setting, adapted to the examination
of wild twistor $\mathcal{D} $-modules, and in much more detail
than needed here, the proof of the above theorem is not so easily
found in loc.cit. It occupies several steps and culminates in the
statement to be found in Remark 7.73, \cite{mochizukiBig}. We plan to
include a short overview of the necessary steps in a future version of this
paper.


\section{The local duality pairing}

Our aim now is to study the meromorphic de Rham complex associated to
the flat algebraic connection $(E, \nabla) $ by lifting to the real
oriented blow-up $\pi: \widetilde{X} \to \Xan $. To this end, we consider
the following de Rham complexes on $\widetilde{X} $:

\begin{definition}
 The {\bf asymptotically flat de Rham complex} is defined to
 be the complex
$$
\DRD(\nabla^\wee) := \AD \otimes_{\pi^{-1}(\O_{\Xan})}
\pi^{-1}(\DR(\nabla^\wee)) \in D^b(\C_{\widetilde{X}}) 
$$
and the {\bf moderate de Rham complex} to be the complex
$$
\DRmod(\nabla) := \Amod \otimes_{\pi^{-1}(\O_{\Xan})}
\pi^{-1}(\DR(\nabla)) \in D^b(\C_{\widetilde{X}}) \ .
$$
\end{definition}

The moderate de Rham complex computes the meromorphic (and hence also
the algebraic) de Rham cohomology of $(E, \nabla) $ (cf. \cite{sabbah1},
Corollaire 1.1.8):
$$
\Hy^k \big( \widetilde{X}, \DRmod( \nabla) \big) = \Hdr^k ( U, E,
\nabla ) \ .
$$

Since multiplication gives a map $\Amod \otimes_\C \AD \to \AD $, we
deduce that the usual wedge product of a differential form with
moderate growth and an asymptotically flat differential form is again
asymptotically flat. This leads to the following

\begin{definition}
The {\bf local duality pairing} is the natural pairing
\begin{equation} \label{eq:locdu} \DRmod( \nabla^\vee ) \otimes_\C
  \DRD( \nabla ) \to \DRD(\O_X, d)
\end{equation}
induced by the wedge product of forms and the natural contraction $E
\otimes E^\vee \to \O_X $ of the vector bundle $E $ and its dual to the
trivial line bundle.
\end{definition}

Before stating the local duality result, we want to prove that the
complexes involved are concentrated in one degree whenever $(E,
\nabla) $ admits a good lattice on $(X,D) $. 

\begin{proposition} \label{claim:deg0} 
If $(E, \nabla) $ admits a good lattice on $(X, D) $, the complexes
$\DRmod( \nabla ) $ and $\DRD( \nabla ) $ have cohomology in degree
zero only. 
\end{proposition}
\begin{proof} 
Due to Mochizuki's results, Theorem \ref{thm:mochiformal} and Theorem
\ref{thm:mochiasymp} above, there exists a bicyclic ramification
$\rho:Y \to X $ such that locally on $\widetilde{Y} $ the pull-back
connection $\rho^{-1}(\nabla) $ admits a decomposition
\eqref{eq:asymplift} for small enough multisectors $S $. Let 
$\pi_X:\widetilde{X} \to \Xan $ and $\pi_Y:\widetilde{Y} \to Y^{\rm an} $ denote
the oriented real blow-up of $Z:=\rho^{-1}(D) $ and $D $
respectively. The projection formula for the lift $\widetilde{\rho}:\widetilde{Y}
\to \widetilde{X} $ gives
\begin{equation} \label{eq:rhoDR}
\bR \widetilde{\rho}_\ast \DRmY (\rho^{-1} \nabla) = \bR
\widetilde{\rho}_\ast \AmdY \Lotimes_{\pi^{-1}(\O_X)} \pi_X^{-1}
\DR(\nabla) \ .
\end{equation}
Now, $\widetilde{\rho} $ being a finite map and since obviously $\bR
\widetilde{\rho}_\ast \AmdY= \Amod $ (using the resolution
\eqref{eq:Amodresol}), it follows that it suffices to
prove the claim on $\widetilde{Y} $. Hence, we can assume 
that we have the decompositions \ref{eq:asymplift} for the connection
$(E, \nabla) $ locally on $\widetilde{X} $. Since the decomposition is
$\nabla $-flat, we can further restrict to the case of one summand
$\mathcal{G}_{\fra, S} $, i.e. we only have to consider the case where
a priori $(E, \nabla) $ is of the form $(E, \nabla)= e^{\fra} \otimes
R_{\fra} $ for some regular singular connection $R_{\fra} $. Since
every regular singular connection is a successive extension of regular
singular line bundles, we can further reduce to the case
$R_\fra=x^\lambda $ with a $\lambda \in \C^d $.

The proof will now follow the same arguments in the asymptotically
flat as well as in the moderate case. Since it is a 
local statement, we restrict to the local
situation at some point $x_0 \in D =\{ x_1 \cdots x_k=0 \} $. Let $\vt
\in \pi^{-1}(x_0) \simeq (S^1)^k $ be a direction in $\hD $ over $x_0
$. Then the complex of stalks at $\vt $ which we have to consider 
is given as
\begin{equation} \label{eq:PDE} \ldots \lto \Big( \Asn{? D_t}
  \otimes_{\pi^{-1}\O_X} \pi^{-1} \Omega_X^p \Big)_\vt
  \stackrel{\nabla}{\lto} \Big( \Asn{? D}
  \otimes_{\pi^{-1}\O_X} \pi^{-1} \Omega_X^{p+1} \Big)_\vt
  \lto \ldots \ ,
\end{equation}
where ? stands for either $< $ or mod. In degree $p $ and with the
usual basis $dx_I $ for $I=\{ 1 \le i_1 < \ldots < i_p \le d\} $ of
$\Omega_X^p $, the connection map $\nabla $ in \eqref{eq:PDE} reads as
$$
\sum_{\# I=p} w_I \, dx_I \mapsto \sum_{\# J=p+1} \Big( \sum_{j \in J}
sgn_J(j) (Q_j w_{J \smallsetminus \{ j \}}) \Big) dx_J
$$
with
$$
Q_j u:= \frac{\partial}{\partial x_j} u + \frac{\partial
  \fra}{\partial x_j} \cdot u + x_j^{-1} \lambda_j \cdot u \ ,
$$
and where we define $sgn_J(j) := (-1)^\nu $ for $J=\{ j_1 <
\ldots < j_{p+1} \} $ and $j_\nu=j $.

Consider a germ of a section $\omega $ of $\Asn{? D_t} \otimes
\pi^{-1} \Omega_{X_t}^{p+1} $,
$$
\omega = \sum_{\# J= p+1} w_J \, dx_J \ ,
$$
such that $\nabla \omega = 0 $. We have to find a $p $-form $\eta $
with appropriate growth condition such that $\nabla \eta= \omega $.

To achieve this, let $s \in \N $ be an integer such that $w_J=0 $ for all
$J $ with $J \cap \{ 1, \ldots, s-1 \} \neq \emptyset $ (which is an
empty condition for $s=1 $). We prove that we can find a $p $-form
$\eta $ with coefficients in $\Asn{? D} $ such that
\begin{equation} \label{eq:induction} 
\big( \omega - \nabla \eta \big) \in \sum_{J \cap \{1, \ldots, s\}=
  \emptyset} \Ass{\vt}{? D} \, dx_J \ .
\end{equation}
The assertion then follows immediately by induction.

By assumption
\begin{equation} \label{eq:omegaK} 0=\nabla \omega = \sum_{\#K=p+2}
  \Big( \sum_{k \in K} sgn_K(k) (Q_k w_{K \smallsetminus \{ k\}})
  \Big) dx_K \ .
\end{equation}
Taking $k<s $ and $s \in J $ and examining the summand of
\eqref{eq:omegaK} corresponding to a set of the form $K:= \{ q \} \cup
J $ we see
that
\begin{equation}\label{eq:intgby}
  Q_q w_J = 0 \mbox{\quad for all such } q=1, \ldots, s-1 \ .
\end{equation}

Now, consider the system $(\Sigma_J) $ of partial differential
equations for the unknown function $u_J $, where $J $ is a fixed
subset $J \subset \{ s, \ldots, d\} $ of cardinality $p+1 $ with $s
\in J $:
$$
(\Sigma_J): \left\{
  \begin{array}{ll}
    Q_q u_J= 0 & \mbox{for all } q=1, \ldots, s-1 \\[0.3cm]
    Q_s u_J= w_J \ ,
  \end{array}\right.
$$
together with the integability assumption \eqref{eq:intgby}. Systems
of this type had been studied by Majima (cf. \cite{majima}) before. In the
case of $\Ass{\vt}{<D} $, the result follows from \cite{sabbah2},
Appendix A, in the moderate case we refer to \cite{hien}, Theorem A.1
(which is formulated in dimension two only but generalizes without
difficulties to the case of arbitrary dimension). 

In each case, we can always find a solution $u_J \in \Ass{\vt}{?  D}
$ for any such $J \subset \{s, \ldots, d\} $ of cardinality $p+1 $ and
$s \in J $ and if we let
$$
\eta:= \sum_{J \text{ as above}} u_J \, dx_J \ ,
$$
we easily see that \eqref{eq:induction} is satisfied.
\end{proof}

In particular, the proposition applies to the trivial conenction
$(\O_x, d) $ for any $(X,D) $ such that $D $ has normal
crossings. The only interesting cohomology sheaf is
$$
\H^0 \, \DRD( \O_X, d) = \widetilde{\jmath}_! \C_{\Uan} \ ,
$$
where $\widetilde{\jmath}: \Uan \hookrightarrow \widetilde{X} $ again
denotes the inclusion. Therefore, the local duality pairing can be
written as
\begin{equation} \label{eq:locdupair_j}
\DRmod( \nabla^\vee ) \otimes_\C \DRD( \nabla ) \to
\widetilde{\jmath}_! \C_{\Uan} \ .
\end{equation}

We are now ready to state the following local duality theorem which
lies in the heart of this work:
\begin{theorem} \label{thm:locduality} 
The local duality pairing \eqref{eq:locdu} is perfect in the sense
that the induced morphisms 
$$
\DRmod( \nabla^\vee ) \to \bR \Hom_{\hX} (\DRD( \nabla ),
\widetilde{\jmath}_! \C) \mbox{\quad and \quad} \DRD( \nabla ) \to \bR
\Hom_{\hX} (\DRmod( \nabla^\vee ), \widetilde{\jmath}_! \C)
$$
are isomorphisms in the derived category.
\end{theorem}
\noindent We give the proof of this theorem in the next section.

\subsection{Proof of Theorem \ref{thm:locduality}}

According to Mochizuki's results, $(E, \nabla) $ admits a good lattice
after some regular birational morphism $b:Y \to X $. Let $Z:=b^{-1}(D)
$ be the resulting divisor in $Y $ and let $\pi_X:\widetilde{X} \to \Xan $ and 
$\pi_Y:\widetilde{Y} \to Y^{\rm an} $ denote the oriented real
blow-ups. Let $\widetilde{b}:\widetilde{Y} \to \widetilde{X} $ be a lift
of $b $. We will also use the following notation for the embeddings
$\widetilde{\jmath}:\Xan \smallsetminus D \hookrightarrow 
\widetilde{X} $ and $\widetilde{\imath}: Y^{\rm an} \smallsetminus Z
\hookrightarrow \widetilde{Y} $. The de Rham complexes on the real
oriented blow-ups behave functorial with respect to this situation:
\begin{lemma} \label{lemma:blow}
  We have natural isomorphisms
$$
\bR\widetilde{b}_\ast (\DRmY(b^\ast \nabla)) \cong \DRmod(\nabla)
\mbox{ and } \bR\widetilde{b}_\ast (\DRdY(b^\ast \nabla)) \cong
\DRD(\nabla) \ .
$$
\end{lemma}
\begin{proof} 
Via the projection formula -- note that $E $ is locally free -- we obtain
$$
\bR \widetilde{b}_\ast ( \dr{Y}{? Z} (b^\ast \nabla) ) = \bR
\widetilde{b}_\ast \big( \As{Y}{? Z} \Lotimes_{\pi_Y^{-1}(\O_{Y^{\rm
      an}})} \pi_Y^{-1} {\rm DR}_{Y^{\rm an}}(b^\ast \nabla) \big)
\cong \bR \widetilde{b}_\ast ( \As{Y}{? Z})
\Lotimes_{\pi_X^{-1}(\O_{\Xan})} \pi_X^{-1} \DR(\nabla) \ ,
$$
where $? $ stands for $< $ or mod respectively. Hence, it suffices to
prove that $\bR \widetilde{b}_\ast ( \As{Y}{? Z} ) = \As{X}{? D} $.

Since the restiction of $b:Y^{\rm an} \to \Xan $ to $Y^{\rm an}
\smallsetminus Z $ is an isomorphism  $Y^{\rm an} \smallsetminus Z \to \Xan
\smallsetminus D= \Uan $, we see that the sheaves of flat $C^\infty
$-functions on each space, i.e. all of whose derivations vanish on the
boundary, are compatible with $\widetilde{b} $, i.e.
$$
\widetilde{b}_\ast {\mathcal P}_{\widetilde{Y}}^{<Z} = {\mathcal
  P}_{\widetilde{X}}^{<D} \ .
$$
Similarily, we have $\widetilde{b}_\ast {\mathcal
  P}_{\widetilde{Y}}^{{\rm mod } Z} = {\mathcal
  P}_{\widetilde{X}}^{{\rm mod } D} $ for the $C^\infty $-functions
with moderate growth for arguments in $\Uan $ approaching the
boundary. Computing $\bR \widetilde{b}_\ast ( \As{Y}{? Z} ) $ via the
flat resolutions \eqref{eq:ADresol} and \eqref{eq:Amodresol} 
respectively, the assertions of the lemma follow immediately.
\end{proof}

Now $\widetilde{\imath}_! \C_{\Uan}[2d] $ is the dualizing sheaf
on the manifold $\widetilde{Y} $ with boundary $\widetilde{Y}
\smallsetminus \Uan $ and similarly $\widetilde{\jmath}_!
\C_{\Uan}[2d] $ for $\widetilde{X} $. By Proposition \ref{claim:deg0}
for the trivial connection $(\O_Y, d) $, we deduce from
\eqref{eq:ADresol} the fine resolution
$$
\widetilde{\imath}_! \C_{\Uan} \simeq
\big( {\mathcal P}_{\widetilde{Y}}^{< Z}
\otimes_{\pi_Y^{-1}(C^\infty_{Y^{\rm an}})} \pi_Y^{-1}
\EYpq{\bdot}{\bdot} , \partial, \dolb \, \big) \ ,
$$
where the right hand is to be understood as the simple complex
associated to the indicated double complex of Dolbeault-type. Applying
$\widetilde{b}_\ast $ thus yields an isomorphism
$$
\alpha:\bR \widetilde{b}_\ast \widetilde{\imath}_! \C_{\Uan} =
\widetilde{b}_\ast \big( {\mathcal P}_{\widetilde{Y}}^{< Z}
\otimes_{\pi_Y^{-1}(C^\infty_{Y^{\rm an}})} \pi_Y^{-1}
\EYpq{\bdot}{\bdot} \big) = \widetilde{b}_\ast {\mathcal
  P}_{\widetilde{Y}}^{< Z} \otimes_{\pi_X^{-1}(C^\infty_{\Xan})}
\pi_X^{-1} \Epq{\bdot}{\bdot} \cong \widetilde{\jmath}_! \C_{\Uan} \ .
$$
Since the isomorphisms of Lemma \ref{lemma:blow} had been constructed
in the same way, the diagram
\begin{equation} \label{eq:blowwedge}
\begin{CD}
\bR \widetilde{b}_\ast \DRmY (b^\ast \nabla) \otimes_\C \bR
\widetilde{b}_\ast \DRdY(b^\ast \nabla^\wee) @>>> \bR \widetilde{b}_\ast
\widetilde{\imath}_! \C_{\Uan} \\
@V{\text{Lemma \ref{lemma:blow}}}V{\cong}V @V{\cong}V{\alpha}V \\
\DRmod(\nabla) \otimes_\C \DRD(\nabla^\wee) @>>> \widetilde{\jmath}_!
\C_{\Uan}
\end{CD}
\end{equation}
commutes. Hence, the morphism $\DRmod(\nabla) \to \bR
\cHom_{\widetilde{X}} \big( \DRD(\nabla^\wee), \widetilde{\jmath}_!
\C_{\Uan} \big) $ induced by the lower row of \eqref{eq:blowwedge}
factors as
\begin{equation} \label{eq:rhomwedge}
\begin{CD}
\bR \widetilde{b}_\ast \DRmY (b^\ast \nabla) @>{\beta}>> \bR
\widetilde{b}_\ast \bR \cHom_{\widetilde{Y}} \big( \DRdY(b^\ast
\nabla^\wee), \widetilde{\imath}_! \C_{\Uan} \big) \\
@V{\cong}VV @VV{\gamma}V \\
\DRmod(\nabla) @>>> \bR \cHom_{\widetilde{X}} \big( \DRD(\nabla^\wee),
\widetilde{\jmath}_! \C_{\Uan} \big) \ ,
\end{CD}
\end{equation}
where $\gamma $ is given by the composition of the natural morphism
$$
\bR \widetilde{b}_\ast \bR \cHom_{\widetilde{Y}} \big( \DRdY,
\widetilde{\imath}_! \C_{\Uan} \big) \to 
\bR \cHom_{\widetilde{X}} \big( \bR \widetilde{b}_\ast \DRdY, \bR
\widetilde{b}_\ast \widetilde{\imath}_! \C_{\Uan} \big)
$$
with the morphism $\alpha $ from above. By Poincar\'e-Verdier duality
(Proposition 3.1.10 in \cite{kascha}), $\gamma $ is an
isomorphism. The same arguments apply in the case of the moderate
de Rham complex. Consequently, in order to prove Theorem
\ref{thm:locduality} namely the fact that the bottom row is an isomorphism, it
suffices to do so for $\beta $. 

In other words, we may a priori assume that $(E,
\nabla) $ has a good lattice on $(X,D) $. Then we know by
Proposition \ref{claim:deg0} that both complexes in the local duality
pairing have cohomology in degree $0 $ only, i.e. we are reduced to
prove that the pairing of sheaves
$$
\Smod \otimes \SD \to \widetilde{\jmath}_! \C_{\Uan}
$$
is perfect, where $\Smod:= \H^0(\DRmod( \nabla^\vee )) $ and
$\SD:=\H^0(\DRD( \nabla )) $. 

Exactly as in the proof of Proposition \ref{claim:deg0}, due to the
existence of a good lattice on $(X,D) $ and the asymptotic
decomposition given by Theorem \ref{thm:mochiasymp}, we can reduce to
the case where $\nabla $ is a rank one connection of the elementary
form $\nabla= e^\alpha \otimes x^\lambda $ for some $\alpha \in \O_X(
\ast D) $ and $\lambda \in \C^d $. 

We restrict to a small enough open polysector $V \subset
\hX \stackrel{\pi}{\to} X $ where locally $D=\{ x_1\cdots x_k=0 \}
$. We will write $\alpha(x)= x_1^{-m_1} \cdot x_k^{-m_k} \cdot u(x) $
with non-vanishing $u(x) $. Let us define the Stokes multi-directions
of $\alpha $ along $D $ inside $V $ to be 
\begin{equation} \label{eq:stokesdir}
\Sigma_{\alpha}:= St^{-1} \big( ( \frac{\pi}{2}, \frac{3 \pi}{2})
\big) \ ,
\end{equation}
where $St: \hD \cap V \to \R/ 2\pi \Z $ denotes the map
$$
St(r_i, \vt_i) := - \sum_{i=1}^k m_i \vt_i + \arg(u \circ \pi(r_i,
\vt_i)) \ .
$$
Let $V_{\alpha}:= (V \smallsetminus \hD) \cup \Sigma_{\alpha} $, i.e.
$V_{\alpha} \cap \hD $ are the directions in which $e^{\alpha(x)} $ has
rapid decay for $x $ radially approaching $D $. If
$\widetilde{\jmath}_{\alpha}:V_{\alpha} \hookrightarrow V $ denotes the
inclusion, one obviously has (possibly after shrinking $V $):
\begin{equation} \label{eq:Valpha} \Smod|_V =
  \widetilde{\jmath}_{-\alpha}(e^{-\alpha(x)} \cdot \C_{\Uan})|_V \mbox{\quad
    and \quad} \SD|_V = \widetilde{\jmath}_{\alpha} ( e^{\alpha(x)} \cdot
  \C_{\Uan})|_V \ .
\end{equation}
Consequently we have the following commutative diagramm:
\begin{equation} \label{eq:locV}
  \begin{CD}
    \bR \cHom ( \SD, \widetilde{\jmath}_!
    \C_{\Uan} )|_V @>>> \Smod|_V  \\
    @V{\cong}VV @VV{\cong}V \\
    \bR \cHom \big( (\widetilde{\jmath}_{-\alpha})_! \C_{V \smallsetminus
      D} \, , \, \widetilde{\jmath}_! \C_{V \cap \Uan} \big)|_V @>>>
    (\widetilde{\jmath}_{\alpha})_! \, \C_{V_{\alpha}}
  \end{CD}
\end{equation}
By the factorization $\widetilde{\jmath} = \widetilde{\jmath}_{-\alpha}
\circ \widetilde{\iota}_{-\alpha} $ with $\widetilde{\iota}_{-\alpha}: V
\cap \Uan \hookrightarrow V_{-\alpha} $, we see that
\begin{multline*}
  \bR \cHom \big( (\widetilde{\jmath}_{-\alpha})_!  \C_{V_{-\alpha}} \, , \,
  \widetilde{\jmath}_! \C_{\Uan} \big) \cong
  (\widetilde{\jmath}_{-\alpha})_\ast \bR \cHom \big( \C_{V_{-\alpha}} ,
  (\widetilde{\iota}_{-\alpha})_! \C_{V \cap \Uan} \big)  \cong \\
  \cong (\widetilde{\jmath}_{-\alpha})_\ast \cHom \big( \C_{V_{-\alpha}} ,
  (\widetilde{\iota}_{-\alpha})_! \C_{V \cap \Uan} \big) =
  (\widetilde{\jmath}_{-\alpha})_! \, \C_{V_{\alpha}} \ ,
\end{multline*}
since $(V \smallsetminus V_{\alpha}) \cap \hD $ coincides with the
closure of $V_{-\alpha} \cap \hD $ inside $\hD $.  Hence, the bottom line
of \eqref{eq:locV} is an isomorphism and thus
$$
\Smod \cong \bR \cHom_{\hX} \big( \SD, \widetilde{\jmath}_!  \C_{\Uan}
\big)
$$
locally on $\hX $ over an arbitrary point of $D $.  Interchanging
$\SD $ and $\Smod $ gives the analogous isomorphism.

\noindent
This completes the proof of the local duality result, Theorem
\ref{thm:locduality}. \\
\qed

\section{Period integrals} \label{sec:periodintegrals}

\subsection{Definition of rapid decay homology} \label{sec:defofrd}

We now want to interpret the local duality pairing as a pairing via
period integrals. To this end, let us recall the definition of rapid
decay homology in \cite{hien}.

Note that for the mere definition in \cite{hien}, we consider the
following geometric set-up. We start with the given flat algebraic
connection $(E, \nabla) $ on the smooth quasi-projective variety $U
$. We then compactify $U $ by some smooth projective $X $ with a
normal crossing divisor $D:=X \smallsetminus U $ as complement. Due to
Mochizuki's Theorem, quoted as Theorem \ref{thm:mochiformal} above, we
can pull-back the connection with respect to some regular birational
map $b:Y \to X $ and obtain a good lattice for the connection. In
other words, by replacing $X $ with $Y $ as a different choice of
compactification, we {\bf may assume that the given connection
  admits a good lattice on the chosen compactification $(X, D) $}.

Let $\pi:\widetilde{X} \to X $ denote
the real oriented blow-up of the normal crossing divisor $D \subset X
$. We will write $\widetilde{D}:= \pi^{-1}(D) $ and denote by $j:U \to
X $ and $\widetilde{\jmath}: U \to \widetilde{X} $ the inclusions. For
any $p \in \N_0 $, we will write $S_p(\widetilde{X}) $ for the free
$\Q $-vector space over all piecewise smooth maps $c: \Delta^p \to
\widetilde{X} $ and we will consider singular homology with $\Q
$-coefficients and piecewise smooth chains in the following. 

Let $\Ctd^{-p} $ denote the sheaf of relative $p $-chains of the pair
$\widetilde{D} \subset \widetilde{X} $, i.e. the sheaf associated to
the presheaf 
$$
V \mapsto S_p \big( \widetilde{X}, ( \widetilde{X}  \smallsetminus V) \cup
\widetilde{D} \big) \ .
$$
Let $\bve := \ker(\nabla|_U) $ be the local system on $U $ of flat
sections of $\nabla $. According to \cite{hien}, Definition 2.3, a
local section $c \otimes \ve \in \Gamma(V, \Ctd^{-\bdot} \otimes_\Q
\widetilde{\jmath}_\ast \bve) $ is called a {\em rapid decay chain},
if for any point $y \in c(\Delta^p) \cap
\widetilde{D} \cap V $ the following holds: Let $x_1, \ldots, x_d $ be
local coordinates of $X $ around $y=0 $ such that $D= \{ x_1 \cdots
x_k=0 \} $ for some $1 \le k \le d $. We chose a meromorphic basis $e:
E \cong (\O_X(\ast D))^r $ of $E $ at $y $ and require that if we
develop $\ve = \sum_{i=1}^d f_i \cdot e_i $
in this basis with analytic functions $f_i \in j_\ast \O_U $, these
coefficients $f_i(x) $ decrease faster than any monomial for $x $
approahcing $D $, i.e. that for
all $N \in \N^k $ there is a $C_N>0 $ such that for all small enough
$x $ we have 
$$
| f_i(x) | \le C_N \cdot |x_1|^{N_1} \cdots |x_k|^{N_k} \ .
$$
For chains $c \otimes \ve $ inside $U $, we do not impose any
condition. The sheaves of all these rapid decay $p $-chains will be denoted
by $\Crdtp{-p}(\nabla) $. Together with the usual boundary operator
$\partial $ of singular chains they define the {\em complex of 
rapid decay chains} $\Crdt(\nabla) $. The {\em rapid decay homology}
of $(E(\ast D), \nabla) $ is the hypercohomology
$$
\Hrd_k(\Uan, E,\nabla) := \Hy^{-k}(\widetilde{X}, \Crdt(\nabla)) \ ,
$$
which can be computed as the cohomology of the global sections,
since the usual barycentric subdivision operator can be defined on the
rapid decay chains and thus one deduces that the complex of rapid
decay chains, similar to the sheaf complex of singular chains
(cp.~\cite{swan}, p.~87) is homotopically fine. For more details, we
refer to \cite{hien}, section 2.3.

\subsection{The local duality pairing and periods}
\label{sec:locduperi}

In order to obtain the desired interpretation of the local duality
pairing \eqref{eq:locdu} in terms of period integrals, we will
examine the relation between the asymptotically flat de Rham
complex of the dual connection $(E^\vee, \nabla^\vee) $ and its rapid
decay complex. Let $\bve^\vee $ denote the local system of the dual
connection.

Recall that we chose $(X,D) $ in a way such that $(E, \nabla) $
admits a good lattice on $(X,D) $. According to Proposition
\ref{claim:deg0}, the natural inclusion 
\begin{equation} \label{eq:DRDSD}
\DRD(\nabla^\vee) \stackrel{\sim}{\longleftarrow}
\H^0(\DRD(\nabla^\vee))=: \SD 
\end{equation}
is a quasi-isomorphism in this situation. Note that by definition and
the properties of $\AD $, the sheaf $\SD $ consists of local solutions
of $\nabla^\vee $ which have rapid decay along $D $. More precisely,
given an open subset $V \subset \hX $ and a section $\sigma \in
\Gamma(V.\SD) $, the section $\sigma \in \SD(V) \subset
\widetilde{\jmath}_\ast (\bve^\vee ) $ is rapidly decaying along any
chain in $c \in \Ctd^{-\bdot}(V) $, i.e.
$$
c \otimes \sigma \in \Gamma(V,\Crdt( \nabla^\vee)) \ .
$$
Thus, we have a natural morphism
\begin{equation} \label{eq:ctopSDCrd}
\Ctd^{-\bdot} \otimes_\Q \SD \lto \Crdt( \nabla^\vee ) \ .
\end{equation}
Now, by standard arguments (e.g. \cite{verdier}), the complex
$\Ctd^{\bdot} $ is a homotopically fine resolution of the
constant sheaf $\C_{\widetilde{X}}[2d] $: to see this, note first that $U:=
\widetilde{X} \smallsetminus \widetilde{D} $ is a real manifold of
dimension $2d $ and since we take homology relative to the boundary
$\widetilde{D} $ we deduce that for small enough $V \subset \widetilde{X} $
such that $V \subset \Omega $ for some contractible $\Omega \subset
\widetilde{X} $ we have by excision:
$$
\Gamma \big( V, \mathcal{H}^{-p}( \Ctd^{\bdot}) \big) = H_p \big( \Omega,
(\Omega \smallsetminus V) \cup \widetilde{D} \big) \cong \left\{
  \begin{array}{ll} 0 & \mbox{for } p \neq 2d \\
\C & \mbox{for } p=2d \ .
\end{array} \right.
$$
Consequently, we end up with a natural morphism
$$
\SD[2d] \stackrel{\simeq}{\lto} \C_{\widetilde{X}}[2d] \otimes \SD
\stackrel{\simeq}{\lto} \Ctd^{-\bdot} \otimes \SD
\stackrel{\eqref{eq:ctopSDCrd}}{\lto} \Crdt (\nabla^\vee) \ .
$$

\begin{proposition} \label{prop:CrdDRD}
If $(E, \nabla) $ admits a good lattice on $(X,D) $, there is a natural isomorphism
$$
\Crd (\nabla^\vee) \cong \DRD( \nabla^\vee) [2d] 
$$
in the derived category $D^b( \C_{\widetilde{X}}) $.
\end{proposition}
\begin{proof}
Due to \eqref{eq:DRDSD}, it remains to prove that \eqref{eq:ctopSDCrd}
is a quasi-isomorphism. Restricted to $U $, this is clear since there
is no rapid decay condition involved and the rapid decay homology
is nothing but the usual singular homology with values in the
local system $\bve^\vee $ which coincides with the restriction $\SD|_U
$ on $U $. 

As for the situation locally around $\widetilde{D} $: By assumption on
the existence of a good lattice, there is a cyclic ramification map $\rho:Y
\to X $ such that $\rho^\ast \nabla^\vee $ admits an unramified good
lattice. Since \eqref{eq:ctopSDCrd} obviously is
compatible with direct image by $\rho $, it is enough to prove the
assertion on $Y $, i.e. we can assume without loss of generality that
$\nabla^\vee $ itself has an unramified good lattice. Exactly as in
the proof of Proposition \ref{claim:deg0}, we can apply Mochizuki's
lifting theorem, Theorem \ref{thm:mochiasymp} above, and reduce the
statement to the case $\nabla^\vee = e^\fra \otimes \mathcal{R} $ for some
regular singular connection $\mathcal{R} $. Now, since the solutions
of a regular singular connection are moderate (see \cite{deligne}), we have
$$
\SDo (e^\fra \otimes \cR ) = \SDo (e^\fra) \otimes
\widetilde{\jmath}_\ast \cR 
$$
with the obvious notation $\SDo( \nabla) := \mathcal{H}^0( \DRD
(\nabla )) $ for a given connection $\nabla $. 

Now, consider a small open polysector $V $ around some direction $\vt
\in \pi^{-1}(y) $ with $y \in D $. Let $\Sigma_{\fra} \subset
\pi^{-1}(y) $ denote the set of Stokes-directions of $\fra $ as in
\eqref{eq:stokesdir}. 

We distinguish the following cases with respect to the direction $\vt
$: If $\vt \in \Sigma_{\fra} $ then for $V $ being a small enough
polysector, we have $\SDo(e^\fra) |_V = \widetilde{\jmath}_! \big(
e^{\fra} \C_U \big) |_V $. For a smooth topological chain $c $ in $\hX $,
the local section $e^{\fra} $ will not have rapid decay along $c $ in
$V $ as required by the definition unless the chain does not meet $\hD
\cap V $. Hence
$$
\Crdt( e^{\fra})|_V = \Ctd^{-\bdot} \otimes \SDo( e^\fra)|_V \ .
$$

If $\vt \in \Sigma_{\fra} $, we can assume that $V $ is an open
polysector such that all the arguments of points in $V $ are contained
in $\Sigma_{\fra} $. Then $\SDo( e^{\fra}) |_V \cong
\widetilde{\jmath}_\ast (e^{\fra} \C_V) $.  Similarly, all twisted
chains $c \otimes e^{\fra} $ will have rapid decay inside $V $ and
again both complexes considered are equal to $\Ctd^{-\bdot} \otimes
\widetilde{\jmath}_\ast (e^{\fra} \C_U) $.

Finally, if $\vt $ separates the Stokes regions of $e^{\fra} $ and
$e^{-\fra} $, we have with the notation of \eqref{eq:Valpha}:
$$
\SDo (e^{\fra}) |_V \cong (\widetilde{\jmath}_{\fra})_!
(e^{\fra} \C_{U})|_V \ . 
$$
The subspace $V_{\fra} $ is characterized by the property that
$V_{\fra} \cap \hD $ consists of those directions along which
$e^{\fra(x)} $ has rapid decay for $x $ approaching $\hD $.  In
particular, $c \otimes e^{\fra} $ is a rapid decay chain on $V $ if
and only if the topological chain $c $ in $\hX $ approaches $\hD \cap
V $ in $V_{\fra} $ at most. Hence
$$
\Crdt(e^{-\fra})|_V = \Ctd^{-\bdot} \otimes ( \widetilde{\jmath}_{\fra}
)_!  (e^{\fra} \C_{U}) = \Ctd^{-\bdot} \otimes \SD|_V \ .
$$
Thus, we have seen that \eqref{eq:ctopSDCrd} is a quasi-isomorphism in
all three cases.
\end{proof}

We remain in the situation given above, namely $(E, \nabla) $
admitting a good lattice on $(X,D) $. Due to Proposition
\ref{claim:deg0}, the local duality pairing reduces to the pairing
\begin{equation} \label{eq:locdureduced}
\SD \otimes \Smod \lto \widetilde{\jmath}_! \C_{U} 
\end{equation}
given by multiplication.

In order to be able to apply Proposition \ref{prop:CrdDRD}, consider
the resolution
$$
\AD \otimes_{\pi^{-1}\O_{X}} \pi^{-1} \Omega_{X}^r
\hookrightarrow \big( \mathcal{P}_{\hX}^{<D} \otimes_{\pi^{-1}
  C_{X}^\infty} \pi^{-1} \Epq{r}{\bdot}, \dolb \big) \ ,
$$
induced by the resolution \eqref{eq:ADresol}, where $\Epq{r}{s} $
denotes the sheaf of $C^\infty $ forms on $X $ of degree $(r,s) $.
Let $\Rd^{\bdot} $ denote the total complex associated to the bicomplex 
$$
\Rd^{\bdot, \bdot}:= \big( \mathcal{P}_{\hX}^{<D}
\otimes_{\pi^{-1}C_{X}^\infty} \pi^{-1} \Epq{\bdot}{\bdot} \, ,
\partial, \dolb \big) \ .
$$
Then, $\Rd^{\bdot} $ is a fine resolution of $\widetilde{\jmath}_!
\C_{U} $. 

Next, consider the canonical quasi-isomorphism $\beta: \Ctd^{-\bdot}
\otimes \Rd^{\bdot} \stackrel{\simeq}{\lto} \Dbrd{-\bdot} $ which maps
an element $c \otimes \rho \in \Ctd^{-r} \otimes \Rd^{s}(V) $ of the
left hand side over some open $V \subset \hX $ to 
the distribution given by $\eta \mapsto \int_c \eta \wedge \rho $ for
a test form $\eta $ with compact support in $V $. Because of the rapid
decay property of $\rho $, this distribution has the property
expressed in \eqref{eq:rddistrib}, hence is indeed a local section of
$\Dbrd{s-r} $.

Now, the main step towards the interpretation as period integrals is
to define the following pairing:
\begin{equation} \label{eq:integral}
\gamma: \Crdtp{-r}(\nabla^\wee)  \otimes \DRmd{s}(\nabla) \to
\Dbrd{s-r} \ , \ 
( c \otimes \ve, \omega) \mapsto \big( \eta \mapsto \int_c \eta \wedge
\langle \ve, \omega \rangle \big) \ ,
\end{equation}
where we observe that since $\ve $ is rapidly decaying along $c $ and
$\omega $ is of moderate growth at most, the integral first of all
converges and it satisfies the estimate \eqref{eq:rddistrib} for all
test forms $\eta $ with compact support on $\widetilde{X} $. The
resulting pairing \eqref{eq:integral} is called the {\em period
  pairing} of the connection $\nabla $.

Now, the inclusion $\Smod \to \DRmod (\nabla) $ is a quasi-isomorphism and the
period pairing and the local duality pairing (in the form of
\eqref{eq:locdureduced} and with the appropriate shift) fit into the
commutative diagramm
\begin{equation}\label{eq:finalcommdiag}
  \begin{CD}
    \SD[2d] \otimes \Smod @>>> \widetilde{\jmath}_! \C_{U}[2d] \\
    @V{\simeq}VV @VV{\simeq}V \\
    \Ctd^{-\bdot} \otimes \SD \otimes \Smod @>>> \Ctd^{-\bdot}
    \otimes \Rd^{\bdot} \\
    @V{\simeq}V{\mbox{\eqref{eq:ctopSDCrd}}}V
    @V{\simeq}V{\beta}V \\
    \Crdt( \nabla^\vee ) \otimes \DRmod ( \nabla ) @>{\gamma}>>
    \Dbrd{-\bdot} \ . 
  \end{CD}
\end{equation}
Applying the local duality result, Theorem \ref{thm:locduality}, together
with Poincar\'e-Verdier duality gives an isomorphism
\begin{equation}\label{eq:PV}
  \bR \Gamma \big(\DRD( \nabla^\vee )[2d] \big) \cong \bR \Gamma \bR
  \cHom_{\C_{\widetilde{X}}} ( \DRmod( \nabla ),
  \widetilde{\jmath}_!\C_{U}[2d]) \cong \Hom_\C^{\bdot} \big(\bR
  \Gamma ( \DRmod( \nabla )), \C \big) \ .
\end{equation}

In summary, we have obtained a description of the perfect duality
$$
\Hy^{2d-p} ( \widetilde{X}, \DRD(\nabla^\vee)) \otimes \Hy^p (
\widetilde{X}, \DRmod ( \nabla)) \lto \C
$$
induced by the local duality pairing, in terms of the period pairing
\eqref{eq:integral}. We summarize this result:

\begin{theorem} \label{thm:periodpair}
  The period pairing \eqref{eq:integral} induces a perfect pairing
$$
\Hrd_p ( U, E^\vee, \nabla^\vee) \otimes_\C \Hdr^p( U; E, \nabla )
\lto \C
$$
betweem the algebraic de Rham cohomology of $(E, \nabla) $ and the
rapid decay homology of the dual connection.
\end{theorem}

\subsection{The period determinant}

Let $k \subsetneq \C $ be a subfield of the field of complex numbers
and assume that the given geometric data, the variety $U $ as well as
the vector bundle and the connection are defined over $k $
already. Then the de Rham cohomology inherits this $k $-structure,
i.e. $\Hdr(U, E, \nabla) $ is naturally a $k $-vector space.  

Now, let $F \subsetneq \C $ denote another subfield of $\C $ and let
us assume that the local system $\bve $ on $U^{an} $ comes equipped
with a given $F $-structure. In analogy to \cite{saitotera} (for the
regular singular case), we
consider the category $W_{k,F}(U) $ -- which we already defined in
\cite{hien}, 2.5 -- of triples $\M=((E,\nabla),\bve_F, \rho) $ with:
\begin{enumerate} 
\item a vector bundle $E $ on $U $ with rank $r $ together with a flat
  connection $\nabla:E \to E \otimes_{\O_U} \Omega^1_U $,
\item a local system $\bve_F $ of $F $-vector spaces on the analytic
  manifold $U^{an} $,
\item a morphism $\rho: \bve_F \to E^{an} $ of sheaves on $U^{an} $
  inducing an isomorphism $\bve_F \otimes_F \C \stackrel{\sim}{\to}
  \ker(\nabla^{an}) $ of local systems of $\C $-vector spaces on
  $U^{an} $.
\end{enumerate}
A morphism between $((E, \nabla), \bve_F,
\rho) $ and $((E', \nabla'), {\bve_F}', \rho') $ is given by a
morphism $E \to E' $ respecting the connections together with a
morphism $\bve_F \to {\bve_F}' $ of $F $-local systems with the
natural compatibility condition with respect to $\rho $ and $\rho' $.

Let $((E,\nabla), \bve_F, \rho) \in W_{k,F}(U) $ be an object in this
category for given subfields $k,F \subsetneq \C $. Then the local
system $\bve^\vee $ inherits an $F $-structure from the given $F
$-structure on $\bve $. Consequently, we can consider the $F $-lattice 
$$
\Crdtp{-p}(\nabla^\wee)_F \subset \Ctd^{-p} \otimes_{\Q}
\widetilde{\jmath}_\ast \bve^\wee_F
$$
of all rapidly decaying chains in $\Ctd^{-p} \otimes_{\Q}
\widetilde{\jmath}_\ast \bve^\wee_F $ and end up with a natural $F
$-lattice $\Hrd_p(\Uan, E, \nabla)_F $ inside the rapid decay homology:
$$
\Hrd_p(\Uan, E, \nabla)_F \otimes_F \C \stackrel{\cong}{\lto} \Hrd_p(\Uan,
E, \nabla) \ ,
$$
the isomorphism induced by $\rho $. The duality between the algebraic
de Rham cohomology and the rapid decay homology via the period pairing
enables us to compare these lattices and to generalize \cite{hien},
Definition 2.7, unconditionally to the case of arbitrary dimension of $U $:
\begin{definition}\label{def:det}
  For $((E, \nabla), \bve_F, \rho) \in W_{k,F}(U) $, we define its
  {\bf period determinant} to be the element
$$
\det((E,\nabla), \bve_F, \rho) := \prod_{p \ge 0} \det ( \langle
\gamma^{(p)}_j , \omega^{(p)}_i \rangle )_{i,j}^{(-1)^p} \in \C^\times
/ k^\times F^\times \ , 
$$
where $\omega^{(p)}_i $ denotes a basis of $\Hdr^p(U, E, \nabla) $
over $k $ and $\gamma^{(p)}_j $ a basis of the $F $-vector space
$\Hrd_k(\Uan, E^\wee, \nabla^\wee)_F $.  
\end{definition}

We hope that this will be the starting point of further examinations
in the spirit of T. Saito and T. Terasoma's work \cite{saitotera} in
the case of regular singular connections now for irreguar ones also.


\begin{thebibliography}{3}

\bibitem{andre} Y.~Andr\'e, {\it Comparison theorems between algebraic
    and analytic De Rham cohomology (with emphasis on the $p $-adic
    case)}, J.~de Th\'eorie des Nombres de Bordeaux 16 (2004), 335 --
  355

\bibitem{bbe} A.~Beilinson, S.~Bloch, H.~Esnault, {\it $\ve $-factors
    for Gauss-Manin determinants}, Mosc.~Math.~J. 2(2002), no.~3, 477
  -- 532
  
\bibitem{b-e} S.~Bloch, H.~Esnault, {\it Homology for irregular
    connections}, Journal de Th\'eorie des Nombres de Bordeaux 16
  (2004), 65 -- 78
  
\bibitem{b-e:irreg} S.~Bloch, H.~Esnault, {\it Gauss-Manin determinant
    connections and periods for irregular connections}, GAFA 2000 (Tel
  Aviv 1999), Geom.~Funct.~Anal.~2000, Special Volume, Part I, 1 -- 31
  
\bibitem{b-e:gm} S.~Bloch, H.~Esnault, {\it Gauss-Manin determinants
    for rank $1 $ irregular connections on curves. With an appendix in
    French by P.~Deligne}, Math.~Ann.~321(2001), no.~1, 15 --87
  
\bibitem{borel} A.~Borel, {\it Algebraic ${\mathcal D} $-modules},
  Perspectives in Math.~2, Academic Press, Boston 1987

\bibitem{bredon} G.E.~Bredon, {\it Sheaf Theory}, Second Edition,
  Grad.~texts in Math.~170, Springer-Verlag, Berlin Heidelberg, 1997

\bibitem{deligne} P.~Deligne, {\it Equations diff\'erentielles \`a
    points singuliers r\'eguliers}, LNS 163, Springer-Verlag, Berlin
  Heidelberg, 1970
  
\bibitem{groth} A.~Grothendieck, {\it On the de Rham cohomology of
    algebraic varieties}, Publ.~Math.~IHES 29 (1966), 93 -- 103
  
\bibitem{hien} M. Hien, {\it Periods for irregular singular
    connections on surfaces}, Math. Ann. 337 (2007), 631 -- 669

\bibitem{hoerm} L. H\"ormander, {\it The analysis of linear partial
    differential operators, I}, Springer-Verlag, Berlin-Heidelberg,
  1990

\bibitem{kascha} M.~Kashiwara, P.~Schapira, {\it Sheaves on
    manifolds}, Grundlehren der mathem.~Wissenschaft 292, Springer
  Verlag Berlin Heidelberg, 1990

\bibitem{majima} H.~Majima, {\it Asymptotic analysis for integrable
    connections with irregular singular points}, LNS 1075,
  Springer-Verlag, Berlin, Heidelberg, 1984
  
\bibitem{mal1} B.~Malgrange, {\it Equations diff\'erentielles \`a
    coefficients polynomiaux}, Prog.~in Math.~96, Birkh{\"a}user,
  Basel, Boston, 1991
  
\bibitem{mal3} B.~Malgrange, {\it Ideals of differentiable functions},
  Oxford University Press, 1966

\bibitem{malgrange_reseau} B.~Malgrange, {\it Connexions m\'eromorphes
    II, le r\'eseau canonique}, Invent. Math. 124 (1996), 367 -- 387

\bibitem{mebkhout} Z.~Mebkhout, {\it Le th\'eor\`eme de comparaison
    entre cohomologies de de Rham d'une vari\'et\'e alg\'ebrique
    complexe et le th\'eor\`eme d'existence de Riemann}, Publ.Math.
  IHES 69(1989), 47-89

\bibitem{mebkh2} Z.~Mebkhout, {\it Le formalisme de six op\'erations
    de Grothendieck pour les ${\mathcal D}_X $-modules coh\'erents},
  Travaux en Cours 35, Hermann, Paris, 1989

\bibitem{mochizuki} T.~Mochizuki, {\it Good formal structure for
    meromorphic flat connections on smooth projective surfaces},
  Preprint, March 2008, arXiv:0803.1346v1 [math.AG]

\bibitem{mochizukiBig} T.~Mochizuki, {\it Wild harmonic bundles and
    wild pure twistor $\mathcal{D} $-modules}, Preprint, March 2008,
  arXiv:0803.1344v1 [math.AG]

\bibitem{narv} L.~Narv\'aez Macarro, {\it The local duality theorem in
    ${\mathcal D} $-module theory}, in El\'ements de la th\'eorie des
  syst\`emes diff\'erentiells g\'eom\'etrique, S\'emin.~Congr.~8,
  Soc.~Math.~France, Paris, 2004

\bibitem{ramsib} J.-P.~Ramis, Y.~Sibuya, {\it Hukuhara domains and
    fundamental existence and uniqueness theorems for asymptotic
    solutions of Gevrey type}, Asymptotic Analysis 2 (1989), no.1, 39
  -- 94

\bibitem{sabbah1} C.~Sabbah, {\it Equations diff\'erentielles \`a
    points singuliers irr\'eguliers et ph\'enom\`ene de Stokes en
    dimension 2}, Ast\'erisque 263, Soc.~Math.~de France, 2000
  
\bibitem{sabbah2} C.~Sabbah, {\it Equations diff\'erentielles \`a
    points singuliers irr\'eguliers en dimension 2}, Ann.~Inst.
  Fourier (Grenoble) 43 (1993), 1619 -- 1688

\bibitem{sabbah3} C.~Sabbah, {\it On the comparison theorem for
    elementary irregular ${\mathcal D} $-modules}, Nagoya Math.~J., 141
  (1996), 107 -- 124
  
\bibitem{saitotera} T.~Saito, T.~Terasoma, {\it Determinant of period
    integrals}, J.~AMS 10 (1997), no.~4, 865 -- 937

\bibitem{swan} R.~G.~Swan, {\it The theory of sheaves}, Chicago
  Lectures in Mathematics, The University of Chicago Press,
  Chicago/London, 1964

\bibitem{terasoma} T.~Terasoma, {\it Confluent hypergeometric
    functions and wild ramification}, J.~of Algebra 185 (1996), 1 --
  18

\bibitem{verdier} J.-L.~Verdier, {\it Classe d'homologie associ\'ee a un
    cycle}, in 'S\'eminaire de g\'eom\'etrie analytique',
  Ast\'erisque 36/37, 1976, 101--151

\end{thebibliography}
\end{document}